\documentclass[11pt]{article}
\usepackage{enumerate}
\usepackage{amssymb,a4wide,latexsym,makeidx,epsfig,fleqn}
\usepackage{amsthm}
\usepackage{amsmath}
\usepackage{enumerate}
\usepackage{graphicx}
\usepackage{float}
\allowdisplaybreaks[4]
\newtheorem{theorem}{Theorem}[section]

\newtheorem{lemma}[theorem]{Lemma}

\begin{document}
\textwidth 150mm \textheight 225mm
\title{The sufficient conditions for $k$-leaf-connected graphs in terms of several topological indices\thanks{Supported by the
National Natural Science Foundation of China (No. 12271439).}}
\author{{Tingyan Ma$^{a,b}$, Ligong Wang$^{a,b,}$\thanks{Corresponding author.
E-mail address: lgwangmath@163.com}, Yang Hu$^{a,b}$}\\
{\small $^a$School of Mathematics and Statistics, Northwestern
Polytechnical University,}\\ {\small  Xi'an, Shaanxi 710129,
P.R. China.}\\
{\small $^b$ Xi'an-Budapest Joint Research Center for Combinatorics, Northwestern
Polytechnical University,}\\
{\small Xi'an, Shaanxi 710129,
P.R. China. }\\
{\small E-mail: matingylw@163.com, lgwangmath@163.com, yhumath@163.com} }
\date{}
\maketitle
\begin{center}
\begin{minipage}{120mm}
\vskip 0.3cm
\begin{center}
{\small {\bf Abstract}}
\end{center}
{\small Let $G=(V(G), E(G))$ be a graph with vertex set $V(G)$ and edge set $E(G)$. For $k\geq2$ and given any subset $S\subseteq|V(G)|$ with $|S|=k$, if a graph $G$ of order $|V(G)|\geq k+1$ always has a spanning tree $T$ such that $S$ is precisely the set of leaves of $T$, then the graph $G$ is a $k$-leaf-connected graph. A graph $G$ is called Hamilton-connected if any two vertices of $G$ are connected by a Hamilton path. Based on the definitions of $k$-leaf-connected and Hamilton-connected, we known that a graph is $2$-leaf-connected if and only if it is Hamilton-connected. During the past decades, there have been many results of sufficient conditions for Hamilton-connected with respect to topological indices. In this paper, we present sufficient conditions for a graph $G$ to be $k$-leaf-connected in terms of the Zagreb index, the reciprocal degree distance or the hyper-Zagreb index. Furthermore, we use the first Zagreb index and hyper-Zagreb index of the complement graph $\overline{G}$ to give sufficient conditions for a graph $G$ to be $k$-leaf-connected.
\vskip 0.1in \noindent {\bf Keywords}: $k$-leaf-connected graphs; Zagreb index; Reciprocal degree distance; Hyper-Zagreb index. \vskip
0.1in \noindent {\bf AMS Subject Classification (2020)}: \ 05C09, 05C35.}
\end{minipage}
\end{center}

\section{Introduction }
\label{sec:ch6-introduction}
Throughout this paper, we only consider simple and undirected graphs.
Let $G=(V(G), E(G))$ be a simple graph with vertex set $V(G)=\{v_{1}, v_{2},
\ldots, v_{n}\}$ and edge set $E(G)=\{e_{1},e_{2}, \ldots, e_{m}\}$. We use $e(G)$ to denote the number of edges of $G$. Let $d_{G}(v_{i})$ denote the degree of $v_{i}\in V(G)$. For convenience, we use $\delta=\delta(G)$ to denote the minimum degree of $G$. Without loss of generality, we assume that $d_{G}(v_{1})\leq d_{G}(v_{2})\leq \cdots \leq d_{G}(v_{n})$ is the degree sequence of $G$. A graph $G$ is $d$-regular if $d_{G}(v_{i})=d$ for each $v_{i}\in V(G)$. The distance between $v_{i}$ and $v_{j}$, denotes by $d_{G}(v_{i}, v_{j})$, is the length of a shortest path from $v_{i}$ to $v_{j}$. For two vertex-disjoint graphs $G$ and $H$, the disjoint union of $G$ and $H$, denoted by $G+H$, is the graph with vertex set $V(G)\cup V(H)$ and edge set $E(G)\cup E(H)$. In particular, let $tG$ be the disjoint union of $t$ copies of a graph $G$. The join graph $G\vee H$ is obtained from $G+H$ by adding all possible edges between $V(G)$ and $V(H)$. We use $K_{n}$ to denote a complete graph of order $n$. For undefined terms and notions one can refer to \cite{Bondy2008,Brouwer2011,Godsil2001}.

Topological indices as molecular descriptors, are extensively studied by many scholars and widely used in many fields. They can be used to measure the energy of molecules. Next, we introduce several topological indices which will be used in this paper.

The first Zagreb index $M_{1}(G)$ and the second Zagreb index $M_{2}(G)$ of a graph $G$ are two famous topological indices, which were first came up with Gutman and Trinajsti\'{c} \cite{Gutamn1972} in 1972. They are defined as
\begin{align*}
M_{1}(G)=\sum_{v_{i}\in V(G)}d_{G}^{2}(v_{i})~\rm{and}~\emph{$M_{2}(G)=\sum_{v_{i}v_{j}\in E(G)}d_{G}(v_{i})d_{G}(v_{j})$}.
\end{align*}
Hua and Zhang $\cite{Hua2012}$ introduced a new graph invariant named reciprocal degree distance $RDD(G)$ or additively weighted Harary index \cite{Alizadeh2022}, which is a degree-weighted version of Harary index. The reciprocal degree distance $RDD(G)$ is defined as vertex-degree-weighted sum of reciprocal distances of a connected graph $G$, that is
$$RDD(G)=\sum_{v_{i}\neq v_{j}}\frac{d_{G}(v_{i})+d_{G}(v_{j})}{d_{G}(v_{i}, v_{j})}.$$
Let $\widehat{D}_{G}(v_{i})=\sum_{v_{j}\in V(G)\setminus \{v_{i}\}}\frac{1}{d_{G}(v_{i}, v_{j})}$. Then $$RDD(G)=\sum_{v_{i}\in V(G)}d_{G}(v_{i})\widehat{D}_{G}(v_{i}).$$
By the definition of $\widehat{D}_{G}(v_{i})$, it can easily verified that
$$\widehat{D}_{G}(v_{i})\leq d_{G}(v_{i})+\frac{1}{2}\left(n-1-d_{G}(v_{i})\right).$$

In 2009, Khalifeh, Yousefi-Azari and Ashrafi \cite{Khalifeh2009} firstly proposed the hyper-Zagreb index of $G$, which was defined as
$$HM(G)=\sum_{v_{i}v_{j}\in E(G)}\Big(d_{G}(v_{i})+d_{G}(v_{j})\Big)^{2}$$

The definition of $k$-leaf-connected graphs was firstly introduced by Gurgel and Wakabayashi $\cite{Gurgel1986}$ in 1986, they also gave some properties about the degree and a sufficient condition on edge number of $k$-leaf-connected graphs. For $k\geq2$ and given any subset $S\subseteq|V(G)|$ with $|S|=k$, if a graph $G$ of order $|V(G)|\geq k+1$ always has a spanning tree $T$ such that $S$ is precisely the set of leaves of $T$, then the graph $G$ is a $k$-leaf-connected graph. A graph is called Hamilton-connected if any two vertices of $G$ are connected by a Hamilton path, where a Hamilton path is a path that pass through all the vertices of a graph. According to the definition of $k$-leaf-connected and Hamilton-connected, a graph is $2$-leaf-connected if and only if it is Hamilton-connected.

The study of sufficient conditions of a connected graph $G$ to be Hamiltion-connected can be extended to $k$-leaf-connected. Zhou, Wang and Lu \cite{Zhou2017,Zhou,Zhou2020} obtained some sufficient conditions for a graph to be Hamilton-connected involving the edge number and the (signless Laplacian) spectral radius. Jia, Liu and Du \cite{Jia2018} characterized some sufficient conditions in terms of the Wiener index for a graph to be Hamilton-connected. Later, Zhou, Wang and Lu \cite{Zhou2019} investigated sufficient conditions for a connected graph to be Hamilton-connected involving the Wiener-type invariants. An \cite{An2022} gave sufficient conditions for a graph to be Hamilton-connected by using the first Zagreb index and the reciprocal degree distance.
\noindent\begin{theorem}\label{th1.1} {\normalfont(\cite{An2022})} Let $G$ be a connected graph of order $n\geq6$. If
$$M_{1}(G)> n^{3}-5n^{2}+12n-6,$$
then $G$ is Hamilton-connected. Moreover, $M_{1}(G)=n^{3}-5n^{2}+12n-6$ if and only if $G\cong K_{2}\vee(K_{1}+K_{n-3})$.
\end{theorem}

\noindent\begin{theorem}\label{th1.2} {\normalfont(\cite{An2022})} Let $G$ be a connected graph of order $n\geq6$. If
$$RDD(G)>\frac{1}{2}n(n+1)(2n-1)-(n-1)(5n-6),$$
then $G$ is Hamilton-connected. Moreover, $RDD(G)=\frac{1}{2}n(n+1)(2n-1)-(n-1)(5n-6)$ if and only if $G\cong K_{2}\vee(K_{1}+K_{n-3})$.
\end{theorem}

Lu and Zhou \cite{Lu2022} showed the sufficient conditions for a graph $G$ and its complement graph $\overline{G}$ to be Hamiltion-connected in terms of the first Zagreb index and hyper-Zagreb index of $G$ or its complement in Theorems \ref{th1.3} and \ref{th1.4}.
\noindent\begin{theorem}\label{th1.3} {\normalfont(\cite{Lu2022})} Let $G$ be a connected graph of order $n\geq13$ with minimum degree $\delta\geq3$. If
$$HM(G)\geq2n^{4}-22n^{3}+126n^{2}-306n+372,$$
then $G$ is Hamilton-connected.
\end{theorem}
\noindent\begin{theorem}\label{th1.4}  {\normalfont(\cite{Lu2022})} Let $G$ be a connected graph of order $n\geq6$ with minimum degree $\delta\geq3$. If
$$HM(\overline{G})\leq(n-2)^{2}(2n-9),$$
then $G$ is Hamilton-connected unless $G\in\{K_{3}\vee3K_{1}, K_{6}\vee6K_{1}, K_{5}\vee5K_{1}, K_{4}\vee4K_{1}\}$.
\end{theorem}
For more researches on sufficient conditions of Hamilton-connected graphs, we can see \cite{Chen2018,Wei2019,Yu2013}. On above studies, the sufficient conditions for a graph to be $k$-leaf-connected have been further studied. Ao, Liu, Yuan and Li $\cite{Ao2022}$ improved sufficient conditions for a graph to be $k$-leaf-connected in terms of the edge number and the (signless Laplacian) spectral radius. Moreover, Ma, Ao, Liu, Wang and Hu $\cite{Ma2022}$ further improved the edge condition and the (signless Laplacian) spectral radius conditions for a graph to be $k$-leaf-connected. Recently, Ao, Liu, Yuan and Yu \cite{Ao2023} study some Wiener-type invariants conditions to guarantee a graph to be $k$-leaf-connected. Motivated by the above results. In this paper, we provided sufficient conditions for a graph to be $k$-leaf-connected involving the $M_{1}(G)$, $M_{2}(G)$, $RDD(G)$ and $HM(G)$. Furthermore, we present sufficient conditions for a graph $G$ to be $k$-leaf-connected in terms of $M_{1}(\overline{G})$ and $HM(\overline{G})$.
\section{Preliminaries}
\label{sec:Preliminaries}
We next present some Lemmas that will be used in our subsequent arguments.
\noindent\begin{lemma}\label{le1} {\normalfont(\cite{Gurgel1986})} Let $G$ be a $k$-leaf-connected graph of order $n$, where $2\leq k\leq n-2$. Then $G$ is $(k+1)$-connected.
\end{lemma}
According to Lemma \ref{le1}, we have $\delta\geq k+1$ for any $k$-leaf-connected graph.
\noindent\begin{lemma}\label{le2} {\normalfont(\cite{Gurgel1986})} Assume that $2\leq k\leq n-1$. If $G$ is a graph such that $d_{G}(v_{i})\geq\frac{n+k-1}{2}$ for every vertex $v_{i}$, then $G$ is $k$-leaf-connected.
\end{lemma}
For a $k$-leaf-connected graph, its degree sequence has the following characterization.
\noindent\begin{lemma}\label{le3} {\normalfont(\cite{Gurgel1986})} Let $k$ and $n$ be positive integers such that $2\leq k\leq n-3$. Let $G$ be a graph with degree sequence $d_{G}(v_{1})\leq d_{G}(v_{2})\leq\cdots\leq d_{G}(v_{n})$. Suppose that there is no integer $i$ with $k\leq i\leq \frac{n+k-2}{2}$ such that $d_{G}(v_{i-k+1})\leq i$ and $d_{G}(v_{n-i})\leq n-i+k-2$. Then $G$ is $k$-leaf-connected.
\end{lemma}
The concept of closure of a graph was used implicitly by Ore \cite{Ore1963},
and formally introduced by Bondy and Chvatal \cite{Bondy1976}.
Fix an integer $l\geq 0$, the $l$-closure of a graph $G$ is the graph obtained from $G$
by successively joining pairs of nonadjacent vertices whose degree sum is at least $l$
until no such pair exists. Denote by $C_{l}(G)$ the $l$-closure of $G.$ Then we have
$$d_{C_{l}(G)}(v_{i})+d_{C_{l}(G)}(v_{j})\leq l-1$$ for any pair of nonadjacent vertices $v_{i}$ and $v_{j}$ of $C_{l}(G).$

\begin{lemma}\label{le4} {\normalfont(\cite{Gurgel1986})}
Let $G$ be a graph and $k$ be an integer with $2\leq k \leq n-1.$
Then $G$ is $k$-leaf-connected
if and only if the $(n+k-1)$-closure $C_{n+k-1}(G)$ of $G$ is $k$-leaf-connected.
\end{lemma}

\begin{lemma}\label{le5} {\normalfont(\cite{Ma2022})}
Let $G$ be a connected graph of order $n\geq k+17$ and minimum degree $\delta\geq k+1$, where $k\geq 2$.
If
$$e(G)\geq {n-3\choose 2}+3k+5,$$
then $G$ is $k$-leaf-connected unless $C_{n+k-1}(G)\in\{K_{k}\vee (K_{n-k-2}+K_{2}), K_{3}\vee (K_{n-5}+2K_{1}),
K_{4}\vee (K_{n-7}+3K_{1})\}$.
\end{lemma}

\section{Several Zagreb index conditions for a graph $G$ to be $k$-leaf-connected}
\noindent\begin{theorem}\label{th1} Let $G$ be a connected graph of order $n\geq k+4$ and minimum degree $\delta\geq k+1$, where $2\leq k\leq n-4$. If
$$M_{1}(G)\geq n^{3}-8n^{2}+(4k+25)n+2k^{2}-4k-24,$$
then $G$ is a $k$-leaf-connected graph unless $G\in\{K_{3}\vee(K_{n-5}+2K_{1}), K_{4}\vee4K_{1}, K_{3}\vee (K_{1, 3}+K_{1}), K_{4}\vee(K_{2}+3K_{1}), K_{5}\vee5K_{1}\}$.
\end{theorem}
\begin{proof}
By the contrary, we suppose that $G$ is not $k$-leaf-connected, where $2\leq k\leq n-4$ and $\delta\geq k+1$. Let  $d_{G}(v_{1})\leq d_{G}(v_{2})\leq\cdots\leq d_{G}(v_{n})$ be the degree sequence of $G$. By Lemma \ref{le3}, there exists a positive integer $i$ with $k\leq i\leq \frac{n+k-2}{2}$ such that $d_{G}(v_{i-k+1})\leq i$ and $d_{G}(v_{n-i})\leq n-i+k-2$. Thus, in the following, we consider $k+1\leq\delta\leq d_{G}(v_{i-k+1})\leq i\leq\frac{n+k-2}{2}$. Then
\begin{align*}
M_{1}(G)=&\sum^{n}_{j=1}d_{G}^{2}(v_{j})\\
=&\sum^{i-k+1}_{j=1}d_{G}^{2}(v_{j})+\sum^{n-i}_{j=i-k+2}d_{G}^{2}(v_{j})+\sum^{n}_{j=n-i+1}d_{G}^{2}(v_{j})
\\\leq&(i-k+1)i^{2}+(n-2i+k-1)(n-i+k-2)^{2}+i(n-1)^{2}
\\=&-i^{3}+(5n+4k-8)i^{2}-(3n^{2}-12n+8kn+4k^{2}-14k+11)i+n^{3}+3kn^{2}-5n^{2}+
\\&3k^{2}n-10kn+8n+k^{3}-5k^{2}+8k-4.
\end{align*}
Let
\begin{align*}
f(x)=&-x^{3}+(5n+4k-8)x^{2}-(3n^{2}-12n+8kn+4k^{2}-14k+11)x,
\end{align*}
where $k+1\leq x\leq \frac{n+k-2}{2}$. Taking the first and second derivative of $f(x)$, we get
\begin{align*}
f^{\prime}(x)=&-3x^{2}+2(5n+4k-8)x-3n^{2}+12n-8kn-4k^{2}+14k-11
\end{align*}
and
\begin{align*}
f^{\prime\prime}(x)=-6x+2(5n+4k-8).
\end{align*}
Note that $k+1\leq x\leq\frac{n+k-2}{2}$, that is $n\geq k+4$, we have
\begin{align*}
f^{\prime\prime}(x)=-6x+2(5n+4k-8)\geq7n+5k-10>0.
\end{align*}
Therefore $f(x)$ is a concave function on $x\in[k+1, \frac{n+k-2}{2}]$. Since $x$ is a positive integer, we have $f(x)\leq max\{f(k+1), f(\lfloor\frac{n+k-2}{2}\rfloor)\}$. By directly calculation, we have
\begin{align*}
f(k+1)=-k^{3}-(3n-7)k^{2}-(3n^{2}-14n+12)k-3n^{2}+17n-20.
\end{align*}
If $n+k-2$ is odd, then
\begin{align*}
f\left(\lfloor\frac{n+k-2}{2}\rfloor\right)=&f\left(\frac{n+k-3}{2}\right)\\
=&-\frac{1}{8}\left(3n^{3}+(19k-17)n^{2}+(25k^{2}-78k+29)n+9k^{3}-49k^{2}+71k-15\right).
\end{align*}
Subtracting $f(k+1)$ from $f(\lfloor\frac{n+k-2}{2}\rfloor)$, we get
\begin{align*}
f(k+1)-f(\lfloor\frac{n+k-2}{2}\rfloor)=\frac{1}{8}(n-k-5)(n-k-7)(3n+k-5)>0~\rm{for}~\emph{$n\geq k+8$}.
\end{align*}
If $n+k-2$ is even, then
\begin{align*}
f\left(\lfloor\frac{n+k-2}{2}\rfloor\right)=&f\left(\frac{n+k-2}{2}\right)\\
=&-\frac{1}{8}\left(3n^{3}+(19k-22)n^{2}+(25k^{2}-76k+48)n+9k^{3}-46k^{2}+72k-32\right).
\end{align*}
Subtracting $f(k+1)$ from $f(\lfloor\frac{n+k-2}{2}\rfloor)$, we get
\begin{equation*}
f(k+1)-f(\lfloor\frac{n+k-2}{2}\rfloor)\\
=\frac{1}{8}(n-k-4)[3n^{2}-(2k+34)n-k^{2}-6k+48]>0~\rm{for}~\emph{$n\geq k+10$},
\end{equation*}
Thus $f(x)\leq f(k+1)$, and thus
\begin{align*}
M_{1}(G)\leq n^{3}-8n^{2}+(4k+25)n+2k^{2}-4k-24.
\end{align*}
In the following, we consider the subcases whether $i=k+1$ or not.\\
\textbf{Case 1:} If $i=k+1$. \\
By the above discussion, we have
$$M_{1}(G)=n^{3}-8n^{2}+(4k+25)n+2k^{2}-4k-24.$$
Then all the inequalities in the proof of Theorem \ref{th1} should be equalities. We have $$d_{G}(v_{1})=d_{G}(v_{2})=k+1,~d_{G}(v_{3})=d_{G}(v_{4})=\cdots=d_{G}(v_{n-k-1})=n-3,$$
$$d_{G}(v_{n-k})=d_{G}(v_{n-k+1})=\cdots=d_{G}(v_{n})=n-1.$$
Hence, $G\cong K_{k+1}\vee(K_{n-k-3}+2K_{1})$, which is $k$-leaf-connected for $k\geq3$. For $k=2$, one can check that $G\cong K_{3}\vee(K_{n-5}+2K_{1})$ is not Hamilton-connected, and hence, $G\cong K_{3}\vee(K_{n-5}+2K_{1})$ is not 2-leaf-connected.\\
\textbf{Case 2:} If $i\neq k+1$, we have $k+2\leq i\leq\frac{n+k-2}{2}$.\\
According to the above discussion, we only need to consider other subcases $n=k+6$, $n=k+7$ and $n=k+8$.\\
\textbf{Subcase 2.1:} $n=k+6$.\\
By Lemma \ref{le3}, we have $d_{G}(v_{3})\leq k+2$, $d_{G}(v_{4})\leq k+2$. Then
\begin{equation}\label{eq1}
d_{G}(v_{1})\leq d_{G}(v_{2})\leq d_{G}(v_{3})\leq d_{G}(v_{4})\leq k+2,~d_{G}(v_{5})\leq d_{G}(v_{6})\leq\cdots\leq d_{G}(v_{k+6})\leq k+5.
\end{equation}
Combining with the condition of Theorem \ref{th1} and Inequality (\ref{eq1}), we have
$$k^3+16k^2+57k+54\leq M_{1}(G)=\sum^{k+6}_{j=1}d_{G}^{2}(v_{j})\leq4(k+2)^{2}+(k+2)(k+5)^{2}=k^3+16k^2+61k+66.$$
Thus
\begin{align*}
d_{G}^{2}(v_{5})+d_{G}^{2}(v_{6})+\cdots+d_{G}^{2}(v_{k+6})=M_{1}(G)-\sum^{4}_{j=1}d_{G}^{2}(v_{j})\geq&k^3+16k^2+57k+54-4(k+2)^{2}\\
=&k^3+12k^2+41k+38.
\end{align*}
Since $k(k+5)^{2}+2(k+4)^{2}=k^3+12k^2+41k+32<k^3+12k^2+41k+38$, then there are at least $k+1$ vertices with degree $k+5$ among $d_{G}(v_{5}), d_{G}(v_{6}), \ldots, d_{G}(v_{k+6})$.
According to $d_{G}^{2}(v_{5})\geq k^3+12k^2+41k+38-(k+1)(k+5)^{2}=k^{2}+6k+13$ and $k^3+16k^2+57k+54\leq M_{1}(G)\leq k^3+16k^2+61k+66$, we can get $d_{G}(v_{5})=k+5$ or $d_{G}(v_{5})=k+4$.\\
If $d_{G}(v_{5})=k+5$, then
$$d_{G}(v_{1})=d_{G}(v_{2})=d_{G}(v_{3})=d_{G}(v_{4})=k+2,~d_{G}(v_{5})=d_{G}(v_{6})=\cdots=d_{G}(v_{k+6})=k+5.$$
Hence, $G\cong K_{k+2}\vee4K_{1}$, which is $k$-leaf-connected for $k\geq3$. For $k=2$, one can check that $G\cong K_{4}\vee4K_{1}$ is not Hamilton-connected, and hence, $G\cong K_{4}\vee4K_{1}$ is not 2-leaf-connected.\\
If $d_{G}(v_{5})=k+4$, then
$$d_{G}(v_{1})=k+1,~d_{G}(v_{2})=d_{G}(v_{3})=d_{G}(v_{4})=k+2,~d_{G}(v_{5})=k+4,$$
$$d_{G}(v_{6})=\cdots=d_{G}(v_{k+6})=k+5.$$
Hence, $G\cong K_{k+1}\vee (K_{1, 3}+K_{1})$, which is $k$-leaf-connected for $k\geq3$. For $k=2$, one can check that $G\cong K_{3}\vee (K_{1, 3}+K_{1})$ is not Hamilton-connected, and hence, $G\cong K_{3}\vee (K_{1, 3}+K_{1})$ is not 2-leaf-connected.\\
\textbf{Subcase 2.2:} $n=k+7$.\\
By directly calculation, we can get that $M_{1}(G)=k^{3}+19k^{2}+84k+102$. Then all the inequalities in the proof of Theorem \ref{th1} should be equalities. We have $$d_{G}(v_{1})=d_{G}(v_{2})=d_{G}(v_{3})=k+2,~d_{G}(v_{4})=d_{G}(v_{5})=k+3,$$
$$d_{G}(v_{6})=d_{G}(v_{7})=\cdots=d_{G}(v_{k+7})=k+6.$$
Hence, $G\cong K_{k+2}\vee(K_{2}+3K_{1})$, which is $k$-leaf-connected for $k\geq3$. For $k=2$, one can check that $G\cong K_{k+2}\vee(K_{2}+3K_{1})$ is not Hamilton-connected, and hence, $G\cong K_{4}\vee(K_{2}+3K_{1})$ is not 2-leaf-connected.\\
\textbf{Subcase 2.3:} $n=k+8$.\\
Note that $k+2\leq i\leq \frac{n+k-2}{2}=k+3$. For $i=k+2$, we have
$$d_{G}(v_{1})\leq d_{G}(v_{2})\leq d_{G}(v_{3})\leq k+2,~d_{G}(v_{4})\leq d_{G}(v_{5})\leq d_{G}(v_{6})\leq k+4,$$
$$d_{G}(v_{7})\leq d_{G}(v_{8})\leq\cdots\leq d_{G}(v_{k+8})\leq k+7,$$
while $3(k+2)^{2}+3(k+4)^{2}+(k+2)(k+7)^{2}=k^{3}+22k^{2}+113k+158<k^{3}+22k^{2}+117k+176$, a contradiction. \\
Then $i=k+3$. By Lemma \ref{le3}, we have $d_{G}(v_{4})\leq k+3$, $d_{G}(v_{5})\leq k+3$. Then
\begin{equation}\label{eq2}
d_{G}(v_{1})\leq d_{G}(v_{2})\leq\cdots\leq d_{G}(v_{5})\leq k+3,~d_{G}(v_{6})\leq d_{G}(v_{7})\leq\cdots\leq d_{G}(v_{k+8})\leq k+7.
\end{equation}
Combining with the condition of Theorem \ref{th1} and Inequality (\ref{eq2}), we have
$$k^{3}+22k^{2}+117k+176\leq M_{1}(G)=\sum^{k+8}_{j=1}d_{G}^{2}(v_{j})\leq5(k+3)^{2}+(k+3)(k+7)^{2}=k^3+22k^{2}+121k+192.$$
Thus
\begin{align*}
d_{G}^{2}(v_{6})+d_{G}^{2}(v_{7})+\cdots+d_{G}^{2}(v_{k+8})=M_{1}(G)-\sum^{5}_{j=1}d_{G}^{2}(v_{j})\geq&k^{3}+22k^{2}+117k+176-5(k+3)^{2}\\
=&k^{3}+17k^{2}+87k+131.
\end{align*}
Since $(k+1)(k+7)^{2}+2(k+6)^{2}=k^{3}+17k^{2}+87k+121<k^{3}+17k^{2}+87k+131$, then there are at least $k+2$ vertices with degree $k+7$ among $d_{G}(v_{6}), d_{G}(v_{7}), \cdots, d_{G}(v_{k+8})$.
According to $d_{G}^{2}(v_{6})\geq k^{3}+17k^{2}+87k+131-(k+2)(k+7)^{2}=k^{2}+10k+33$ and $k^{3}+22k^{2}+117k+176\leq M_{1}(G)\leq k^3+22k^{2}+121k+192$, we can $d_{G}(v_{6})=k+7$ and
$$d_{G}(v_{1})=d_{G}(v_{2})=\cdots=d_{G}(v_{5})=k+3,~d_{G}(v_{6})=d_{G}(v_{7})=\cdots=d_{G}(v_{k+8})=k+7.$$
Hence, $G\cong K_{k+3}\vee5K_{1}$, which is $k$-leaf-connected for $k\geq3$. For $k=2$, one can check that $G\cong K_{5}\vee5K_{1}$ is not Hamilton-connected, and hence, $G\cong K_{5}\vee5K_{1}$ is not 2-leaf-connected.\\
\end{proof}

\noindent\begin{theorem}\label{th2} Let $G$ be a connected graph of order $n\geq k+4$ and minimum degree $\delta\geq k+1$, where $2\leq k\leq n-4$. If
$$RDD(G)\geq n^{3}-7n^{2}+(20+4k)n+k^{2}-4k-17,$$
then $G$ is a $k$-leaf-connected graph unless $G\in\{K_{3}\vee(K_{n-5}+2K_{1}), K_{4}\vee4K_{1}, K_{3}\vee (K_{1, 3}+K_{1}), K_{4}\vee(K_{2}+3K_{1}), K_{5}\vee5K_{1}\}$.
\end{theorem}
\begin{proof}
By the contrary, we suppose that $G$ is not $k$-leaf-connected, where $2\leq k\leq n-4$ and $\delta\geq k+1$. Let  $d_{G}(v_{1})\leq d_{G}(v_{2})\leq\cdots\leq d_{G}(v_{n})$ be the degree sequence of $G$. By Lemma \ref{le3}, there exists a positive integer $i$ with $k\leq i\leq \frac{n+k-2}{2}$ such that $d_{G}(v_{i-k+1})\leq i$ and $d_{G}(v_{n-i})\leq n-i+k-2$. Thus, in the following, we consider $k+1\leq\delta\leq d_{G}(v_{i-k+1})\leq i\leq\frac{n+k-2}{2}$. Then
\begin{align*}
RDD(G)=&\sum_{v_{i}\in V(G)}d_{G}(v_{i})\widehat{D}_{G}(v_{i})\\
\leq&\sum_{v_{i}\in V(G)}d_{G}(v_{i})\left(\frac{1}{2}\left(n-1-d_{G}(v_{i})\right)\right)\\
=&\frac{1}{2}\left((n-1)\sum_{v_{i}\in V(G)}d_{G}(v_{i})+\sum_{v_{i}\in V(G)}d_{G}^{2}(v_{i})\right)\\
\leq&\frac{1}{2}(n-1)[(i-k+1)i+(n-2i+k-1)(n-i+k-2)+i(n-1)]+\\
&\frac{1}{2}[(i-k+1)i^{2}+(n-2i+k-1)(n-i+k-2)^{2}+i(n-1)^{2}]\\
=&\frac{1}{2}[-i^{3}+(8n+4k-11)i^{2}-(5n^{2}-19n+12kn+4k^{2}-18k+16)i+\\
&2n^{3}+(5k-9)n^{2}+(4k^{2}-15k+13)n+k^{3}-6k^{2}+11k-6].
\end{align*}
Let
\begin{align*}
g(x)=&-x^{3}+(8n+4k-11)x^{2}-(5n^{2}-19n+12kn+4k^{2}-18k+16)x,
\end{align*}
where $k+1\leq x\leq \frac{n+k-2}{2}$. Taking the first and second derivative of $g(x)$, we get
\begin{align*}
g^{\prime}(x)=&-3x^{2}+2(8n+4k-11)x-(5n^{2}-19n+12kn+4k^{2}-18k+16)
\end{align*}
and
\begin{align*}
g^{\prime\prime}(x)=-6x+2(8n+4k-11).
\end{align*}
Note that $x\leq\frac{n+k-2}{2}$, we have
\begin{align*}
g^{\prime\prime}(x)=-6x+2(8n+4k-11)\geq13n+5k-16>0~~(n\geq k+4).
\end{align*}
Therefore $g(x)$ is a concave function on $x\in[k+1, \frac{n+k-2}{2}]$. Since $x$ is a positive integer, we have $g(x)\leq max\{g(k+1), g(\lfloor\frac{n+k-2}{2}\rfloor)\}$. By directly calculation, we have
\begin{align*}
g(k+1)=-k^{3}-(4n-8)k^{2}-(5n^{2}-23n+19)k-5n^{2}+27n-28.
\end{align*}
If $n+k-2$ is odd, then
\begin{align*}
g\left(\lfloor\frac{n+k-2}{2}\rfloor\right)=&g\left(\frac{n+k-3}{2}\right)\\
=&-\frac{1}{8}\left(5n^{3}+(31k-27)n^{2}+(35k^{2}-122k+43)n+9k^{3}-59k^{2}+103k-21\right).
\end{align*}
Subtracting $f(k+1)$ from $f(\lfloor\frac{n+k-2}{2}\rfloor)$, we get
\begin{align*}
f(k+1)-f(\lfloor\frac{n+k-2}{2}\rfloor)=\frac{1}{8}(n-k-5)(n-k-7)(5n+k-7)>0~\rm{for}~\emph{$n\geq k+8$}.
\end{align*}
If $n+k-2$ is even, then
\begin{align*}
g\left(\lfloor\frac{n+k-2}{2}\rfloor\right)=&g\left(\frac{n+k-2}{2}\right)\\
=&-\frac{1}{8}\left(5n^{3}+(31k-36)n^{2}+(35k^{2}-116k+76)n+9k^{3}-56k^{2}+100k-48\right).
\end{align*}
Subtracting $g(k+1)$ from $g(\lfloor\frac{n+k-2}{2}\rfloor)$ gives
\begin{equation*}
f(k+1)-f(\lfloor\frac{n+k-2}{2}\rfloor)\\
=\frac{1}{8}(n-k-4)[5n^{2}-(4k+56)n-k^{2}-4k+68]>0~\rm{for}~\emph{$n\geq k+10$},
\end{equation*}
which implies that $g(x)\leq g(k+1)$. Thus
\begin{align*}
RDD(G)\leq n^{3}-7n^{2}+(20+4k)n+k^{2}-4k-17.
\end{align*}
In the following, we consider the subcases whether $i=k+1$ or not.\\
\textbf{Case 1:} If $i=k+1$. \\
By the above discussion, we have
$$RDD(G)= n^{3}-7n^{2}+(20+4k)n+k^{2}-4k-17.$$
Then
$$d_{G}(v_{1})=d_{G}(v_{2})=k+1,~d_{G}(v_{3})=d_{G}(v_{4})=\cdots=d_{G}(v_{n-k-1})=n-3,$$
$$d_{G}(v_{n-k})=d_{G}(v_{n-k+1})=\cdots=d_{G}(v_{n})=n-1.$$
Hence, $G\cong K_{k+1}\vee(K_{n-k-3}+2K_{1})$, which is $k$-leaf-connected for $k\geq3$. For $k=2$, one can check that $G\cong K_{3}\vee(K_{n-5}+2K_{1})$ is not Hamilton-connected, and hence, $G\cong K_{3}\vee(K_{n-5}+2K_{1})$ is not 2-leaf-connected.\\
\textbf{Case 2:} If $i\neq k+1$, we have $k+2\leq i\leq\frac{n+k-2}{2}$.\\
According to the above discussion, we only need to consider other subcases $n=k+6$, $n=k+7$ and $n=k+8$.\\
\textbf{Subcase 2.1:} $n=k+6$.\\
By Lemma \ref{le4}, we have $d_{G}(v_{3})\leq k+2$, $d_{G}(v_{4})\leq k+2$. Then
\begin{equation}\label{eq3}
d_{G}(v_{1})\leq d_{G}(v_{2})\leq d_{G}(v_{3})\leq d_{G}(v_{4})\leq k+2,~d_{G}(v_{5})\leq d_{G}(v_{6})\leq\cdots\leq d_{G}(v_{k+6})\leq k+5.
\end{equation}
Combining with the condition of Theorem \ref{th2} and Inequality (\ref{eq3}), we have
\begin{align*}
k^{3}+16k^{2}+64k+67&\leq RDD(G)\\
&\leq\frac{1}{2}\bigg((n-1)\sum_{v_{i}\in V(G)}d_{G}(v_{i})+\sum_{v_{i}\in V(G)}d_{G}^{2}(v_{i})\bigg)\\
&\leq\frac{1}{2}[(n-1)(4(k+2)+(k+2)(k+5))+4(k+2)^{2}+(k+2)(k+5)^{2}]\\
&=k^{3}+16k^2+67k+78.
\end{align*}
Thus
\begin{align*}
&RDD(G)-\frac{1}{2}\bigg((n-1)\sum^{4}_{j=1}d_{G}(v_{j})+\sum^{4}_{j=1}d_{G}^{2}(v_{j})\bigg)\\
\geq &k^{3}+16k^{2}+64k+67-\frac{1}{2}[4(n-1)(k+2)+4(k+2)^{2}]\\
=&k^{3}+12k^{2}+42k+39.
\end{align*}
Since $\frac{1}{2}[(k+5)(k(k+5)+2(k+4))+k(k+5)^{2}+2(k+4)^{2}]=k^{3}+12k^{2}+42k+36<k^{3}+12k^{2}+42k+39$, then there are at least $k+1$ vertices with degree $k+5$ among $d_{G}(v_{5}), d_{G}(v_{6}), \ldots, d_{G}(v_{k+6})$.
According to $\frac{1}{2}[(n-1)d_{G}(v_{5})+d_{G}^{2}(v_{5})]\geq k^{3}+12k^{2}+42k+39-\frac{1}{2}[(n-1)(k+1)(k+5)+(k+1)(k+5)^{2}]=k^{2}+7k+14$ and $k^{3}+16k^{2}+64k+67\leq RDD(G)\leq k^{3}+16k^2+67k+78$, we can get $d_{G}(v_{5})=k+5$ or $d_{G}(v_{5})=k+4$.\\
If $d_{G}(v_{5})=k+5$, then
$$d_{G}(v_{1})=d_{G}(v_{2})=d_{G}(v_{3})=d_{G}(v_{4})=k+2,~d_{G}(v_{5})=d_{G}(v_{6})=\cdots=d_{G}(v_{k+6})=k+5.$$
Hence, $G\cong K_{k+2}\vee4K_{1}$, which is $k$-leaf-connected for $k\geq3$. For $k=2$, one can check that $G\cong K_{4}\vee4K_{1}$ is not Hamilton-connected, and hence, $G\cong K_{4}\vee4K_{1}$ is not 2-leaf-connected.\\
If $d_{G}(v_{5})=k+4$, then
$$d_{G}(v_{1})=k+1,~d_{G}(v_{2})=d_{G}(v_{3})=d_{G}(v_{4})=k+2,~d_{G}(v_{5})=k+4,$$
$$d_{G}(v_{6})=\cdots=d_{G}(v_{k+6})=k+5.$$
Hence, $G\cong K_{k+1}\vee (K_{1, 3}+K_{1})$, which is $k$-leaf-connected for $k\geq3$. For $k=2$, one can check that $G\cong K_{3}\vee (K_{1, 3}+K_{1})$ is not Hamilton-connected, and hence, $G\cong K_{3}\vee (K_{1, 3}+K_{1})$ is not 2-leaf-connected.\\
\textbf{Subcase 2.2:} $n=k+7$.\\
By directly calculation, we can get that $RDD(G)=k^{3}+19k^{2}+93k+123$. We have $$d_{G}(v_{1})=d_{G}(v_{2})=d_{G}(v_{3})=k+2,~d_{G}(v_{4})=d_{G}(v_{5})=k+3,$$
$$d_{G}(v_{6})=d_{G}(v_{7})=\cdots=d_{G}(v_{k+7})=k+6.$$
Hence, $G\cong K_{k+2}\vee(K_{2}+3K_{1})$, which is $k$-leaf-connected for $k\geq3$. For $k=2$, one can check that $G\cong K_{k+2}\vee(K_{2}+3K_{1})$ is not Hamilton-connected, and hence, $G\cong K_{4}\vee(K_{2}+3K_{1})$ is not 2-leaf-connected.\\
\textbf{Subcase 2.3:} $n=k+8$.\\
Note that $k+2\leq i\leq \frac{n+k-2}{2}=k+3$. For $i=k+2$, we have
$$d_{G}(v_{1})\leq d_{G}(v_{2})\leq d_{G}(v_{3})\leq k+2,~d_{G}(v_{4})\leq d_{G}(v_{5})\leq d_{G}(v_{6})\leq k+4,$$
$$d_{G}(v_{7})\leq d_{G}(v_{8})\leq\cdots\leq d_{G}(v_{k+8})\leq k+7,$$
while $\frac{1}{2}[(n-1)(3(k+2)+3(k+4)+(k+2)(k+7))+3(k+2)^{2}+3(k+4)^{2}+(k+2)(k+7)^{2}]=k^{3}+22k^{2}+125k+191<k^{3}+22k^{2}+128k+207$, a contradiction. \\
Then $i=k+3$. By Lemma \ref{le4}, we have $d_{G}(v_{4})\leq k+3$, $d_{G}(v_{5})\leq k+3$. Then
\begin{equation}\label{eq4}
d_{G}(v_{1})\leq d_{G}(v_{2})\leq\cdots\leq d_{G}(v_{5})\leq k+3,~d_{G}(v_{6})\leq d_{G}(v_{7})\leq\cdots\leq d_{G}(v_{k+8})\leq k+7.
\end{equation}
Combining with the condition of Theorem \ref{th2} and Inequality (\ref{eq4}), we have
\begin{align*}
k^{3}+22k^{2}+128k+207&\leq RDD(G)\\
&\leq\frac{1}{2}\bigg((n-1)\sum_{v_{i}\in V(G)}d_{G}(v_{i})+\sum_{v_{i}\in V(G)}d_{G}^{2}(v_{i})\bigg)\\
&\leq\frac{1}{2}[(n-1)(5(k+3)+(k+3)(k+7))+5(k+3)^{2}+(k+3)(k+7)^{2}]\\
&=k^3+22k^{2}+131k+222.
\end{align*}
Thus
\begin{align*}
&RDD(G)-\frac{1}{2}\bigg((n-1)\sum^{5}_{j=1}d_{G}(v_{j})+\sum^{5}_{j=1}d_{G}^{2}(v_{j})\bigg)\\
\geq&k^{3}+22k^{2}+128k+207-\frac{1}{2}[5(n-1)(k+3)+5(k+3)^{2}]\\
=&k^{3}+17k^{2}+88k+132.
\end{align*}
Since $\frac{1}{2}[(n-1)((k+1)(k+7)+2(k+6))+(k+1)(k+7)^{2}+2(k+6)^{2}]=k^{3}+17k^{2}+88k+127<k^{3}+17k^{2}+88k+132$, then there are at least $k+2$ vertices with degree $k+7$ among $d_{G}(v_{6}), d_{G}(v_{7}), \cdots, d_{G}(v_{k+8})$.
According to $\frac{1}{2}[(n-1)d_{G}(v_{6})+d_{G}^{2}(v_{6})]\geq k^{3}+17k^{2}+88k+132-\frac{1}{2}[(n-1)(k+2)(k+7)+(k+2)(k+7)^{2}]=k^{2}+11k+34$ and $k^{3}+22k^{2}+128k+207\leq RDD(G)\leq k^3+22k^{2}+131k+222$, we can get $d_{G}(v_{6})=k+7$ and
$$d_{G}(v_{1})=d_{G}(v_{2})=\cdots=d_{G}(v_{5})=k+3,~d_{G}(v_{6})=d_{G}(v_{7})=\cdots=d_{G}(v_{k+8})=k+7.$$
Hence, $G\cong K_{k+3}\vee5K_{1}$, which is $k$-leaf-connected for $k\geq3$. For $k=2$, one can check that $G\cong K_{5}\vee5K_{1}$ is not Hamilton-connected, and hence, $G\cong K_{5}\vee5K_{1}$ is not 2-leaf-connected.\\
\end{proof}

\noindent\begin{theorem}\label{th3} Let $G$ be a connected graph of order $n\geq k+5$ and minimum degree $\delta\geq k+1$, where $2\leq k\leq n-5$. If
$$HM(G)\geq2n^{4}-22n^{3}+(12k+102)n^{2}-(48k+210)n+4k^{3}+12k^{2}+64k+164,$$
then $G$ is a $k$-leaf-connected graph unless $G\in\{K_{5}\vee5K_{1}, K_{4}\vee4K_{1}\}$.
\end{theorem}
\begin{proof}
By the contrary, we suppose that $G$ is not $k$-leaf-connected, where $2\leq k\leq n-5$ and $\delta\geq k+1$. Let  $d_{G}(v_{1})\leq d_{G}(v_{2})\leq\cdots\leq d_{G}(v_{n})$ be the degree sequence of $G$. By Lemma \ref{le3}, there exists a positive integer $i$ with $k\leq i\leq \frac{n+k-2}{2}$ such that $d_{G}(v_{i-k+1})\leq i$ and $d_{G}(v_{n-i})\leq n-i+k-2$. Thus, in the following, we consider $k+1\leq\delta\leq d_{G}(v_{i-k+1})\leq i\leq\frac{n+k-2}{2}$. Then
\begin{align*}
HM(G)=&\sum_{v_{i}v_{j}\in E(G)}\Big(d_{G}(v_{i})+d_{G}(v_{j})\Big)^{2}\\
\leq&2\sum_{j=1}^{n}d_{G}^{3}(v_{j})\\
\leq&2\left(\sum^{i-k+1}_{j=1}d_{G}^{3}(v_{j})+\sum^{n-i}_{j=i-k+2}d_{G}^{3}(v_{j})+\sum^{n}_{j=n-i+1}d_{G}^{3}(v_{j})\right)
\\\leq&2[(i-k+1)i^{3}+(n-2i+k-1)(n-i+k-2)^{3}+i(n-1)^{3}]
\\=&6i^{4}-(14n+16k-28)i^{3}+(18n^{2}+36kn-66n+18k^{2}-66k+60)i^{2}-\\
&(8n^{3}+30kn^{2}-48n^{2}+30k^{2}n-108kn+90n+10k^{3}-54k^{2}+96k-54)i+\\
&2n^{4}+8kn^{3}-14n^{3}+12k^{2}n^{2}-42kn^{2}+36n^{2}+8k^{3}n-42k^{2}n+72kn-40n+\\
&2k^{4}-14k^{3}+36k^{2}-40k+16.
\end{align*}
Let
\begin{align*}
l(x)=&6x^{4}-(14n+16k-28)x^{3}+(18n^{2}+36kn-66n+18k^{2}-66k+60)x^{2}-\\
&(8n^{3}+30kn^{2}-48n^{2}+30k^{2}n-108kn+90n+10k^{3}-54k^{2}+96k-54)x,
\end{align*}
where $k+1\leq x\leq \frac{n+k-2}{2}$. Taking the first and second derivative of $l(x)$, we get
\begin{align*}
l^{\prime}(x)=&24x^{3}-3(14n+16k-28)x^{2}+2(18n^{2}+36kn-66n+18k^{2}-66k+60)x-\\
&(8n^{3}+30kn^{2}-48n^{2}+30k^{2}n-108kn+90n+10k^{3}-54k^{2}+96k-54)
\end{align*}
and
\begin{align*}
l^{\prime\prime}(x)=72x^{2}-6(14n+16k-28)x+2(18n^{2}+36kn-66n+18k^{2}-66k+60).
\end{align*}
Since $\frac{n+k-2}{2}<-\frac{-6(14n+16k-28)}{2\times72}=\frac{7n+8k-14}{12}$, and
\begin{align*}
l^{\prime\prime}\left(\frac{n+k-2}{2}\right)=12n^{2}+(18k-36)n+6k^{2}-24k+24>0.
\end{align*}
Therefore $l(x)$ is a concave function on $x\in[k+1, \frac{n+k-2}{2}]$. Since $x$ is an integer, we have $l(x)\leq max\{l(k+1), l(\lfloor\frac{n+k-2}{2}\rfloor)\}$. By directly calculation, we have
\begin{align*}
l(k+1)=&-(8k+8)n^{3}-(12k^{2}-54k-66)n^{2}-(8k^{3}-42k^{2}+120k+170)n-\\
&2k^{4}+18k^{3}-24k^{2}+104k+148.
\end{align*}
If $n+k-2$ is odd, then
\begin{align*}
l\left(\lfloor\frac{n+k-2}{2}\rfloor\right)=&l\left(\frac{n+k-3}{2}\right)\\
=&-\frac{7}{8}n^{4}-(\frac{27}{4}k-\frac{29}{4})n^{3}-(12k^{2}-39k+21)n^{2}-(\frac{33}{4}k^{3}-\frac{177}{4}k^{2}+\frac{279}{4}k-\\
&\frac{99}{4})n-\frac{17}{8}k^{4}+\frac{31}{2}k^{3}-\frac{159}{4}k^{2}+\frac{81}{2}k-\frac{81}{8}
\end{align*}
and for $n\geq k+9$
\begin{align*}
&l(k+1)-l\left(\frac{n+k-3}{2}\right)\\
=&\frac{1}{8}(n-k-5)(7n^{3}-(3k+87)n^{2}-(3k^{2}-18k-216)n-k^{3}-15k^{2}-51k-253)>0.
\end{align*}
If $n+k-2$ is even, then
\begin{align*}
l\left(\lfloor\frac{n+k-2}{2}\rfloor\right)=&l\left(\frac{n+k-2}{2}\right)\\
=&-\frac{7}{8}n^{4}-(\frac{27}{4}k-\frac{17}{2})n^{3}-(12k^{2}-21k+27)n^{2}-(\frac{33}{4}k^{3}-\frac{87}{2}k^{2}+72k-\\
&35)n-\frac{17}{8}k^{4}+15k^{3}-39k^{2}+43k-16
\end{align*}
and for $n\geq k+11$
\begin{align*}
&l(k+1)-l\left(\frac{n+k-2}{2}\right)\\
=&\frac{1}{8}(n-k-4)(7n^{3}-(3k+104)n^{2}-(3k^{2}-4k-328)n-k^{3}-20k^{2}-40k-328)>0,
\end{align*}
which implies that $l(x)\leq l(k+1)$. Thus
\begin{align*}
HM(G)\leq2n^{4}-22n^{3}+(12k+102)n^{2}-(48k+210)n+4k^{3}+12k^{2}+64k+164.
\end{align*}
In the following, we consider the subcases whether $i=k+1$ or not.\\
\textbf{Case 1:} If $i=k+1$. \\
By the above discussion, we have
$$HM(G)=2n^{4}-22n^{3}+(12k+102)n^{2}-(48k+210)n+4k^{3}+12k^{2}+64k+164.$$
Then all the inequalities in the proof of Theorem \ref{th3} should be equalities. We have $$d_{G}(v_{1})=d_{G}(v_{2})=k+1,~d_{G}(v_{3})=d_{G}(v_{4})=\cdots=d_{G}(v_{n-k-1})=n-3,$$
$$d_{G}(v_{n-k})=d_{G}(v_{n-k+1})=\cdots=d_{G}(v_{n})=n-1.$$
Hence, $G\cong K_{k+1}\vee(K_{n-k-3}+2K_{1})$, which is $k$-leaf-connected for $k\geq3$. For $k=2$, one can check that $G\cong K_{3}\vee(K_{n-5}+2K_{1})$ is not Hamilton-connected, while it contradicts with $2\sum_{v_{i}v_{j}\in E(G)}d_{G}(v_{i})d_{G}(v_{j})=\sum_{v_{i}v_{j}\in E(G)}(d_{G}^{2}(v_{i})+d_{G}^{2}(v_{j}))$.\\
\textbf{Case 2:} If $i\neq k+1$, we have $k+2\leq i\leq\frac{n+k-2}{2}$.\\
According to the above discussion, we only need to consider other subcases $n=k+6$, $n=k+7$, $n=k+8$ and $n=k+10$.\\
\textbf{Subcase 2.1:} $n=k+6$.\\
By Lemma \ref{le4}, we have $d_{G}(v_{3})\leq k+2$, $d_{G}(v_{4})\leq k+2$. Then
\begin{equation}\label{eq5}
d_{G}(v_{1})\leq d_{G}(v_{2})\leq d_{G}(v_{3})\leq d_{G}(v_{4})\leq k+2,~d_{G}(v_{5})\leq d_{G}(v_{6})\leq\cdots\leq d_{G}(v_{k+6})\leq k+5.
\end{equation}
Combining with the condition of Theorem \ref{th3} and Inequality (\ref{eq5}), we have
\begin{align*}
2k^{4}+42k^{3}+246k^{2}+574k+416\leq HM(G)&\leq2\sum^{k+6}_{j=1}d_{G}^{3}(v_{j})\\
&\leq2(4(k+2)^{3}+(k+2)(k+5)^{3})\\
&=2k^{4}+42k^{3}+258k^{2}+646k+564.
\end{align*}
Thus
\begin{align*}
2(d_{G}^{3}(v_{5})+d_{G}^{3}(v_{6})+\cdots+d_{G}^{3}(v_{k+6}))&=HM(G)-2\sum^{4}_{j=1}d_{G}^{3}(v_{j})\\
&\geq2k^{4}+42k^{3}+246k^{2}+574k+416-8(k+2)^{3}\\
&=2k^{4}+34k^{3}+198k^{2}+478k+352.
\end{align*}
Since $2(k(k+5)^{3}+2(k+4)^{3})=2k^{4}+34k^{3}+198k^{2}+442k+256<2k^{4}+34k^{3}+198k^{2}+478k+352$, then there are at least $k+1$ vertices with degree $k+5$ among $d_{G}(v_{5}), d_{G}(v_{6}), \ldots, d_{G}(v_{k+6})$.
According to $2d_{G}^{3}(v_{5})\geq 2k^{4}+34k^{3}+198k^{2}+478k+352-2(k+1)(k+5)^{3}=2k^{3}+18k^{2}+78k+102$ and $2k^{4}+42k^{3}+246k^{2}+574k+416\leq HM(G)\leq2k^{4}+42k^{3}+258k^{2}+646k+564$, we can get $d_{G}(v_{5})=k+5$ and
$$d_{G}(v_{1})=d_{G}(v_{2})=d_{G}(v_{3})=d_{G}(v_{4})=k+2,~d_{G}(v_{5})=d_{G}(v_{6})=\cdots=d_{G}(v_{k+6})=k+5.$$
Hence, $G\cong K_{k+2}\vee4K_{1}$, which is $k$-leaf-connected for $k\geq3$. For $k=2$, one can check that $G\cong K_{4}\vee4K_{1}$ is not Hamilton-connected, and hence, $G\cong K_{4}\vee4K_{1}$ is not 2-leaf-connected.\\
\textbf{Subcase 2.2:} $n=k+7$.\\
Note that $i=k+2$, we have $d_{G}(v_{3})\leq k+2$, $d_{G}(v_{5})\leq k+3$ and degree sequence
$$d_{G}(v_{1})\leq d_{G}(v_{2})\leq d_{G}(v_{3})\leq k+2,~d_{G}(v_{4})\leq d_{G}(v_{5})\leq k+3,$$
$$d_{G}(v_{6})\leq d_{G}(v_{7})\leq\cdots\leq d_{G}(v_{k+7})\leq k+6.$$
If the degree sequence of $G$ is
$$d_{G}(v_{1})=d_{G}(v_{2})=d_{G}(v_{3})=k+2,~d_{G}(v_{4})=d_{G}(v_{5})=k+3,$$
$$d_{G}(v_{6})=d_{G}(v_{7})=\cdots=d_{G}(v_{k+7})=k+6.$$
We have $G\cong K_{k+2}\vee(K_{2}+3K_{1})$ and
$$HM(G)\leq HM(K_{k+4}\vee6K_{1})=\sum_{v_{i}v_{j}\in E(G)}\Big(d_{G}(v_{i})+d_{G}(v_{j})\Big)^{2}=2k^{4}+50k^{3}+360k^{2}+912k+756,$$\\
which contradicts $HM(G)\geq2k^{4}+50k^{3}+360k^{2}+1044k+948$.\\
\textbf{Subcase 2.3:} $n=k+8$.\\
Note that $k+2\leq i\leq \frac{n+k-2}{2}=k+3$. For $i=k+2$, we have
$$d_{G}(v_{1})\leq d_{G}(v_{2})\leq d_{G}(v_{3})\leq k+2,~d_{G}(v_{5})\leq d_{G}(v_{5})\leq d_{G}({v_{6}})\leq k+4,$$
$$d_{G}(v_{7})\leq d_{G}(v_{8})\leq\cdots\leq d_{G}(v_{k+8})\leq k+7,$$
while $2(3(k+2)^{3}+3(k+4)^{3}+(k+2)(k+7)^{3})=2k^{4}+58k^{3}+486k^{2}+1634k+1804<2k^{4}+58k^{3}+498k^{2}+1742k+1940$, a contradiction. \\
Then $i=k+3$. By Lemma \ref{le4}, we have $d_{G}(v_{4})\leq k+3$, $d_{G}(v_{5})\leq k+3$. Then
\begin{equation}\label{eq6}
d_{G}(v_{1})\leq d_{G}(v_{2})\leq\cdots\leq d_{G}(v_{5})\leq k+3,~d_{G}(v_{6})\leq d_{G}(v_{7})\leq\cdots\leq d_{G}(v_{k+8})\leq k+7.
\end{equation}
Combining with the condition of Theorem \ref{th3} and Inequality (\ref{eq6}), we have
\begin{align*}
2k^{4}+58k^{3}+498k^{2}+1742k+1940\leq HM(G)&\leq2\sum^{k+8}_{j=1}d_{G}^{3}(v_{j})\\
&\leq2(5(k+3)^{3}+(k+3)(k+7)^{3})\\
&=2k^{4}+58k^{3}+510k^{2}+1838k+2328.
\end{align*}
Thus
\begin{align*}
2(d_{G}^{3}(v_{6})+d_{G}^{3}(v_{7})+\cdots+d_{G}^{3}(v_{k+8}))=&HM(G)-\sum^{5}_{j=1}d_{G}^{3}(v_{j})\\
\geq&2k^{4}+58k^{3}+498k^{2}+1742k+1940-10(k+3)^{3}\\
=&2k^{4}+48k^{3}+408k^{2}+1472k+1670.
\end{align*}
Since $2((k+1)(k+7)^{3}+2(k+6)^{3})=2k^{4}+48k^{3}+408k^{2}+1412k+1550<2k^{4}+48k^{3}+408k^{2}+1472k+1670$, then there are at least $k+2$ vertices with degree $k+7$ among $d_{G}(v_{6}), d_{G}(v_{7}), \ldots, d_{G}(v_{k+8})$.
According to $2d_{G}^{3}(v_{6})\geq 2k^{4}+48k^{3}+408k^{2}+1472k+1670-2(k+2)(k+7)^{3}=2k^{3}+30k^{2}+198k+98$ and $2k^{4}+58k^{3}+498k^{2}+1742k+1940\leq HM(G)\leq 2k^{4}+58k^{3}+510k^{2}+1838k+2328$, we can get $d_{G}(v_{6})=k+7$ and
$$d_{G}(v_{1})=d_{G}(v_{2})=\cdots=d_{G}(v_{5})=k+3,~d_{G}(v_{6})=d_{G}({v_{7}})=\cdots=d_{G}(v_{k+8})=k+7.$$
Hence, $G\cong K_{k+3}\vee5K_{1}$, which is $k$-leaf-connected for $k\geq3$. For $k=2$, one can check that $G\cong K_{5}\vee5K_{1}$ is not Hamilton-connected, and hence, $G\cong K_{5}\vee5K_{1}$ is not 2-leaf-connected.\\
\textbf{Subcase 2.3:} $n=k+10$.\\
Note that $k+2\leq i\leq \frac{n+k-2}{2}=k+4$. For $i=k+2$, we have
$$d_{G}(v_{1})\leq d_{G}(v_{2})\leq d_{G}(v_{3})\leq k+2,~d_{G}(v_{4})\leq d_{G}(v_{5})\leq\cdots\leq d_{G}({v_{8}})\leq k+6,$$
$$d_{G}(v_{9})\leq d_{G}(v_{10})\leq\cdots\leq d_{G}(v_{k+10})\leq k+9,$$
while $2(3(k+2)^{3}+5(k+6)^{3}+(k+2)(k+9)^{3})=2k^{4}+74k^{3}+810k^{2}+3582k+5124<2k^{4}+74k^{3}+846k^{2}+4014k+6264$, a contradiction. \\
For $i=k+3$. By Lemma \ref{le4}, we have $d_{G}(v_{4})\leq k+3$, $d_{G}(v_{7})\leq k+5$. Then
$$d_{G}(v_{1})\leq d_{G}(v_{2})\leq d_{G}(v_{3})\leq d_{G}(v_{4})\leq k+3,~d_{G}(v_{5})\leq d_{G}(v_{6})\leq d_{G}({7})\leq k+5,$$
$$d_{G}(v_{8})\leq d_{G}(v_{9})\cdots\leq d_{G}(v_{k+10})\leq k+9,$$
while $2(4(k+3)^{3}+3(k+5)^{3}+(k+3)(k+9)^{3})=2k^{4}+74k^{3}+810k^{2}+3582k+5340<2k^{4}+74k^{3}+846k^{2}+4014k+6264$, a contradiction. \\
For $i=k+4$. By Lemma \ref{le4}, we have $d_{G}(v_{5})\leq k+4$, $d_{G}(v_{6})\leq k+4$. Then
$$d_{G}(v_{1})\leq d_{G}(v_{2})\leq\cdots\leq d_{G}(v_{6})\leq k+4,~d_{G}({v_{7}})\leq d_{G}(v_{8})\leq\cdots\leq d_{G}(v_{k+10})\leq k+9,$$
If the degree sequence of $G$ is
$$d_{G}(v_{1})=d_{G}(v_{2})=\cdots=d_{G}(v_{6})=k+4,~d_{G}({v_{7}})=d_{G}(v_{8})=\cdots=d_{G}(v_{k+10})=k+9.$$
We have $G\cong K_{k+4}\vee6K_{1}$ and
$$HM(G)\leq HM(K_{k+4}\vee6K_{1})=\sum_{v_{i}v_{j}\in E(G)}\Big(d_{G}(v_{i})+d_{G}(v_{j})\Big)^{2}=2k^{4}+74k^{3}+846k^{2}+3678k+5400,$$\\
which contradicts $HM(G)\geq2k^{4}+74k^{3}+846k^{2}+4014k+6264$.
\end{proof}

\noindent\begin{theorem}\label{th4} Let $G$ be a connected graph of order $n\geq k+5$ and minimum degree $\delta\geq k+1$, where $2\leq k\leq n-5$. If
$$M_{2}(G)\geq\frac{1}{2}n^{4}-\frac{11}{2}n^{3}+\frac{6k+51}{2}n^{2}-\frac{24k+105}{2}n+k^{3}+3k^{2}+16k+41,$$
then $G$ is a $k$-leaf-connected graph unless $G\in\{K_{5}\vee5K_{1}, K_{4}\vee4K_{1}\}$.
\end{theorem}
\begin{proof}
We suppose that $G$ is not $k$-leaf-connected, combining the proof of Theorem \ref{th3} with the definition of the second Zagreb index $M_{2}(G)$, we have
\begin{align*}
M_{2}(G)=&\sum_{v_{i}v_{j}\in E(G)}d_{G}(v_{i})d_{G}(v_{j})\\
\leq&\frac{1}{4}\sum_{v_{i}v_{j}\in E(G)}(d_{G}(v_{i})+d_{G}(v_{j}))^{2}\\
\leq&\frac{1}{2}\sum_{i=1}^{n}d_{G}^{3}(v_{i})\\
\leq&\frac{1}{2}n^{4}-\frac{11}{2}n^{3}+\frac{6k+51}{2}n^{2}-\frac{24k+105}{2}n+k^{3}+3k^{2}+16k+41,
\end{align*}
unless $G\in\{K_{5}\vee5K_{1}, K_{4}\vee4K_{1}\}$.
\end{proof}

\section{The first Zagreb index and hyper-Zagreb index conditions on $\overline{G}$ for $G$ to be $k$-leaf-connected}
\noindent\begin{theorem}\label{th5} Let $G$ be a connected graph of order $n\geq k+17$ and minimum degree $\delta\geq k+1$, where $2\leq k\leq n-17$. If
$$M_{1}(\overline{G})\leq(n-k)[3(n-k)-11],$$
then $G$ is $k$-leaf-connected unless $cl_{n+k-1}(G)\in\{K_{k}\vee (K_{n-k-2}+K_{2}), K_{3}\vee (K_{n-5}+2K_{1})\}$.
\end{theorem}
\begin{proof}
By the contrary, we suppose that $G$ is not $k$-leaf-connected, where $2\leq k\leq n-17$ and $\delta\geq k+1$. Let $H=cl_{n+k-1}(G)$, then by Lemma \ref{le4}, $H$ is not $k$-leaf-connected. For any two nonadjacent vertices $v_{i}$ and $v_{j}$ in $H$, we have $d_{H}(v_{i})+d_{H}(v_{j})\leq n+k-2$. Hence, for any two adjacent vertices $v_{i}$ and $v_{j}$ in $\overline{H}$, we have
$$d_{\overline{H}}(v_{i})+d_{\overline{H}}(v_{j})=n-1-d_{H}(v_{i})+n-1-d_{H}(v_{j})\geq n-k.$$
According to the definition of the first Zagreb index, we get
\begin{align*}
M_{1}(\overline{H})=\sum_{v_{i}\in V(\overline{H})}d_{\overline{H}}^{2}(v_{i})=\sum_{v_{i}v_{j}\in E(\overline{H})}(d_{\overline{H}}(v_{i})+d_{\overline{H}}(v_{j}))\geq(n-k)e(\overline{H}).
\end{align*}
Since $H$ is not $k$-leaf-connected graph, by Lemma \ref{le5}, we have $e(H)\leq {n-3\choose 2}+3k+4$ or $e(H)\geq {n-3\choose 2}+3k+5$ and $H\in\{ K_{k}\vee (K_{n-k-2}+K_{2}), K_{3}\vee (K_{n-5}+2K_{1}),
K_{4}\vee (K_{n-7}+3K_{1})\}$.\\
If $e(H)\leq {n-3\choose 2}+3k+4$, note that $\overline{H}\subseteq \overline{G}$, then
\begin{align*}
M_{1}(\overline{G})\geq& M_{1}(\overline{H})\geq(n-k)e(\overline{H})\geq(n-k)\left[{n\choose 2}-{n-3\choose 2}-3k-4\right]\\
=&(n-k)[3(n-k)-10].
\end{align*}
It is easy to check that $(n-k)[3(n-k)-11]<(n-k)[3(n-k)-10]$ for $n\neq k$, this contradicts with the condition of Theorem \ref{th5}.\\
If $e(H)\geq {n-3\choose 2}+3k+5$, then $H\in\{K_{k}\vee (K_{n-k-2}+K_{2}), K_{3}\vee (K_{n-5}+2K_{1}),
K_{4}\vee (K_{n-7}+3K_{1})\}$. By directly calculation, we have
$$M_{1}(\overline{K_{k}\vee (K_{n-k-2}+K_{2})})=2(n-k)(n-k-2),$$
$$M_{1}(\overline{K_{3}\vee (K_{n-5}+2K_{1})})=2(n^{2}-6n+6)$$
and
$$M_{1}(\overline{G})\geq M_{1}(\overline{K_{4}\vee (K_{n-7}+3K_{1})})=3(n^{2}-7n+4).$$
We cannot compare complete $M_{1}(\overline{G})$ with $(n-k)[3(n-k)-11]$ for $H\in\{K_{k}\vee(K_{n-k-2}+K_{2}), K_{3}\vee (K_{n-5}+2K_{1})\}$. Therefore,
$cl_{n+k-1}(G)=H\in\{K_{k}\vee(K_{n-k-2}+K_{2}), K_{3}\vee(K_{n-5}+2K_{1})\}$.
\end{proof}

\noindent\begin{theorem}\label{th6} Let $G$ be a connected graph of order $n\geq k+17$ and minimum degree $\delta\geq k+1$, where $2\leq k\leq n-17$. If
$$HM(\overline{G})\leq(n-k)^{2}[3(n-k)-11],$$
then $G$ is $k$-leaf-connected unless $cl_{n+k-1}(G)\in\{K_{k}\vee (K_{n-k-2}+K_{2}), K_{3}\vee (K_{n-5}+2K_{1})\}$.
\end{theorem}
\begin{proof}
By the contrary, we suppose that $G$ is not $k$-leaf-connected, where $2\leq k\leq n-17$ and $\delta\geq k+1$. Let $H=cl_{n+k-1}(G)$, then by Lemma \ref{le4}, $H$ is not $k$-leaf-connected. For any two nonadjacent vertices $v_{i}$ and $v_{j}$ in $H$, we have $d_{H}(v_{i})+d_{H}(v_{j})\leq n+k-2$. Hence, for any two adjacent vertices $v_{i}$ and $v_{j}$ in $\overline{H}$, we have
$$d_{\overline{H}}(v_{i})+d_{\overline{H}}(v_{j})=n-1-d_{H}(v_{i})+n-1-d_{H}(v_{j})\geq n-k.$$
According to the definition of the hyper-Zagreb index, we get
\begin{align*}
HM(\overline{H})=&\sum_{v_{i}v_{j}\in E(\overline{H})}(d_{\overline{H}}(v_{i})+d_{\overline{H}}(v_{j}))^{2}\geq(n-k)^{2}e(\overline{H}).
\end{align*}
Since $H$ is not $k$-leaf-connected graph, by Lemma \ref{le5}, we have $e(H)\leq {n-3\choose 2}+3k+4$ or $e(H)\geq {n-3\choose 2}+3k+5$ and $H\in\{ K_{k}\vee (K_{n-k-2}+K_{2}), K_{3}\vee (K_{n-5}+2K_{1}),
K_{4}\vee (K_{n-7}+3K_{1})\}$.\\
If $e(H)\leq {n-3\choose 2}+3k+4$, note that $\overline{H}\subseteq \overline{G}$, then
\begin{align*}
HM(\overline{G})\geq& HM(\overline{H})\geq(n-k)^{2}e(\overline{H})\geq(n-k)^{2}\left[{n\choose 2}-{n-3\choose 2}-3k-4\right]\\
=&(n-k)^{2}[3(n-k)-10].
\end{align*}
It is easy to check that $(n-k)^{2}[3(n-k)-10]>(n-k)^{2}[3(n-k)-11]$ for $n\neq k$, this contradicts with the condition of Theorem \ref{th6}.\\
If $e(H)\geq {n-3\choose 2}+3k+5$, then $H\in\{K_{k}\vee (K_{n-k-2}+K_{2}), K_{3}\vee (K_{n-5}+2K_{1}),
K_{4}\vee (K_{n-7}+3K_{1})\}$. Note that $\overline{H}\subseteq \overline{G}$, by directly calculation, we have
$$HM(\overline{K_{k}\vee (K_{n-k-2}+K_{2})})=2(n-k)(n-k-2),$$
$$HM(\overline{K_{3}\vee (K_{n-5}+2K_{1})})=2n^{3}-14n^{2}+16n+24$$
and
$$HM(\overline{G})\geq M_{1}(\overline{K_{4}\vee (K_{n-7}+3K_{1})})=3n^{3}-21n^{2}-24n+216.$$
Hence, they all contradict the condition of Theorem \ref{th6}. Assume that $H\in\{ K_{k}\vee (K_{n-k-2}+K_{2}), K_{3}\vee (K_{n-5}+2K_{1})\}$, we cannot compare complete $HM(\overline{G})$ with $(n-k)^{2}[3(n-k)-11]$. Therefore,
$cl_{n+k-1}(G)=H\in\{K_{k}\vee (K_{n-k-2}+K_{2}), K_{3}\vee (K_{n-5}+2K_{1})\}$.
\end{proof}

\end{document}